\documentclass[11pt,sloppy]{volume}

\usepackage{graphicx}

\usepackage{latexsym}
\usepackage{amsfonts}
\usepackage{amsthm}
\usepackage{amsmath}

\newcommand{\RR}{\mathbb{R}}

\def\operator#1{{\hbox{\bf #1}}}
\def\dom{{\operator{dom}}}
\def\RR{{\operator{R}}}
\def\dmax{{\operator{dmax}}}
\def\dmin{{\operator{dmin}}}

\newtheorem{theorem}{Theorem}[section]
\newtheorem{lemma}[theorem]{Lemma}

\newtheorem{problem}[theorem]{Problem}
\newtheorem{example}[theorem]{Example}

\theoremstyle{definition}

\begin{document}
\title[Some Algorithms in the Kepler Conjecture]{Some algorithms arising
in the proof of the Kepler conjecture}


\toctitle{Some algorithms arising in the proof of the Kepler conjecture}


\author{Thomas C. Hales}

\tocauthor{Thomas C. Hales}
\runningauthor{Thomas C. Hales}

\maketitle

\begin{abstract}
By any account, the 1998 proof of the Kepler conjecture is complex.
The thesis underlying this article is that the proof
is complex because it is highly under-automated.  Throughout that proof,
manual procedures are used where automated ones would have been better
suited.  This paper gives a series of nonlinear optimization algorithms 
and shows how a systematic application of these algorithms would bring
substantial simplifications to the original proof.
\end{abstract}

\section{Introduction}

In 1998 a proof of the Kepler conjecture was completed \cite{H}.  
By any account,
that solution is complex (300 pages
of text, 3GB stored data on the computer,
computer calculations taking months, 40K lines of computer code, and
so forth).  The thesis underlying this article is that the 1998 proof
is complex because it is highly under-automated.  Throughout that proof,
manual procedures are used where automated ones would have been better
suited.
\footnote{Version November 5, 2001}

Ultimately, a properly automated
proof of the Kepler conjecture might be short and elegant.  
The hope is that the Kepler conjecture might eventually become 
an instance of a general
family of optimization problems for which general optimization techniques
exist.  Just as today linear programming problems of a moderate
size can be solved without fanfare, we might hope that problems 
of a moderate size in this
family might be routinely solved by general algorithms.
The proof of the Kepler conjecture would then consist of demonstrating
that the Kepler conjecture can be structured as a problem in this family,
and then invoking the general algorithm to solve the problem.

As a step toward that objective, this article frames the primary algorithms
of that proof in sufficient generality that they may be applied to
much larger families of problems.  The algorithms are arranged into
four sections: 
{\it Quantifier Elimination, Linear Assembly Problems,
 Automated Inequality Proving,
Plane Graph Generation}.

We do not claim any originality in the algorithms.  In fact, the
purpose is just the contrary: to exhibit the proof of the
Kepler conjecture
insofar as possible as an instance of standard optimization techniques.
To keep things as general as possible,
 the algorithms we present here will make no mention
of the particulars of the Kepler conjecture.
A final section lists parts of the 1998 proof 
that can be structured according to these general algorithms.

\section{Quantifier Elimination}
\label{q.e.}

We might try to structure the Kepler conjecture as a
statement in the elementary theory of the real numbers.  Tarski proved
the decidability of this theory, through quantifier
elimination.  G. Collins and others have
developed and implemented concrete algorithms to perform 
quantifier elimination \cite{CJ}.
The Kepler conjecture, as formulated in \cite{FH},
is not a statement in this theory, because the transcendental
arctangent function enters into the statement.  

However, it seems that the arctangent is not essential to the
formulation of the Kepler conjecture, and that it enters only because
no attempt was made to do without it.  For example, it is plausible that
it can be replaced with a
close rational approximation, without doing violence to the
proof.  In fact, the computer calculations in that proof
are already based on rational approximations with explicit error bounds, and
on its rational derivative $1/(1+x^2)$. 

Assuming this can be done, quantifier elimination gives a procedure
to solve the Kepler conjecture.  Unfortunately, these algorithms
are 
prohibitively 
slow (exponential, or doubly exponential in the number of variables).

Section~\ref{linear} of this article proposes a 
different family of optimization problems
for which algorithmic performance is satisfactory.  These are called
linear assembly problems.  

Although quantifier elimination is too slow to be of practical value
as a 1-step solution to the Kepler conjecture, it can be of great value
in proving intermediate results.  Recent algorithms are able to solve
problems nearly at the level of difficulty of intermediate results
in the Kepler conjecture \cite{BPR}, \cite{M}.

Instances of the following families
of problems arise as intermediate steps in the Kepler conjecture.  
Each instance of the following families must provide
an explicit set of parameters $r_k$ (or $\dmin$, $\dmax$), 
and the problem becomes to show that
configurations of points in ${\RR^3}$ with the given parameters do not exist.
In theory, the problems are all amenable to solution by quantifier elimination.

\begin{problem}\label{first} Let $S$ be a simplex whose edges all
have length at most given parameter values $r_i$.  
Show that there is no point in the interior of the
simplex with distance at least $r$ from every vertex.
\end{problem}

\begin{problem} Show that there does not exist a triangle 
of circumradius at most $r_1$, 
and a segment of length at most $r_2$ such that the segment passes through
the interior of the triangle and such that each endpoint of the segment
has distance at least $r_3$ from each vertex of the triangle.
\end{problem}

\begin{problem}\label{pr:third} 
Show that there does not exist a configuration of
$5$ points $0,p_1,p_2,p_3,q$ with given
minimum $\dmin(p,q)$ and maximum $\dmax(p,q)$
distances between each pair $(p,q)$ of points.
The line through $(0,q)$ must link the triangle with vertices $p_i$.
\end{problem}

\begin{problem}
\label{six}
Show that there does not exist a configuration of 
$6$ points $0$, $p_1,\ldots,p_4$, $q$, with given minimum and
maximum distances between each pair $(p,q)$ of points.
The line $(0,q)$ must link the skew quadrilateral with vertices $q_i$
(ordered according to subscripts).
\end{problem}

\begin{problem}
\label{seven}
 Show that there does not exist a configuration of 
$7$ points $0$, $p_1,\ldots,p_4$, $q_1$, $q_2$ 
 with given minimum $\dmin(p,q)$ and
maximum $\dmax(p,q)$ distances between each pair $(p,q)$ of points.
For $j=1,2$, 
the line $(0,q_j)$ must link the skew quadrilateral with vertices $q_i$
(ordered according to subscripts).
\end{problem}

Although we hope that one day these problems will all be amenable
to direct solution by quantifier elimination, in practice, 
we did not try to apply quantifier elimination directly without preprocessing
them.  The idea of preprocessing is that if a configuration exists, then
the points can be moved in such a way to make various upper and lower
bound constraints on distances bind.  With a large number of
binding constraints, the dimension of the configuration space becomes
smaller and the problem easier to solve.  
In some cases, preprocessing reduces the
configuration space to a single configuration, so that the existence of
the configuration can be tested by choosing coordinates and calculating
whether all the metric constraints are satisfied.

These five families of quantifier elimination
 problems have a similar feel to them.
Let us give a preprocessing algorithm in general enough terms that it
applies uniformly to all five problem families.

The primary preprocessing of the configurations is a deformation that we
call {\it pivoting}.  
Fix any three of the points
$p_1$, $p_2$, and $q$ of the configuration.
We call a {\it pivot through axis} 
$(p_1,p_2)$
the continuous 
motion of $q$ in the perpendicular bisecting plane $B$ of
$(p_1,p_2)$ at constant distance from $p_1$ and $p_2$.  Thus, the pivot
moves $q$ in a circular path in the plane $B$.  

\begin{figure}[htb]
  \centering
  \includegraphics{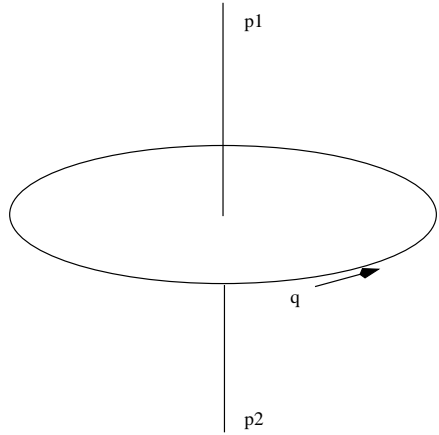}
  \caption{A pivot is the circular motion of a point around a fixed axis.}
  \label{pivot}
\end{figure}

Usually, the plane $P=(p_1,p_2,q)$ 
through the three points is chosen to have the property that
the entire configuration
lies in a half-space through $P$.  
If $q$ moves away from
the half-space containing $P$, the distances from $q$ to the other points
of the configuration increase or remain the same.  
If $q$ moves into the half-space, the
distances decrease or remain the same.

More generally, we allow the plane $P$ to separate the points of the
configuration into two groups, such that the lower distance bounds from
$q$ to the first group do not bind, and such that the upper distance
bounds from $q$ to the second group do not bind.

To apply pivots, we must prescribe their directions, whether into the
half-space or away from it.  To do this, we give a {\it model\/},
which is a configuration that exists, of the form indicated in the
problem family, but which is not required to satisfy the various constraints.  
Various edges in the model are marked with a strut (indicating a lower
bound) or a cable (indicating an upper bound).
If the model
has a cable, then preprocessing pivots are 
applied
to increase the corresponding distance in the configuration space, until
the given upper bound is reached.  Where the model has a strut,
 pivots are applied so as to decrease distances to
the lower bound.

\begin{figure}[htb]
  \centering
  \includegraphics{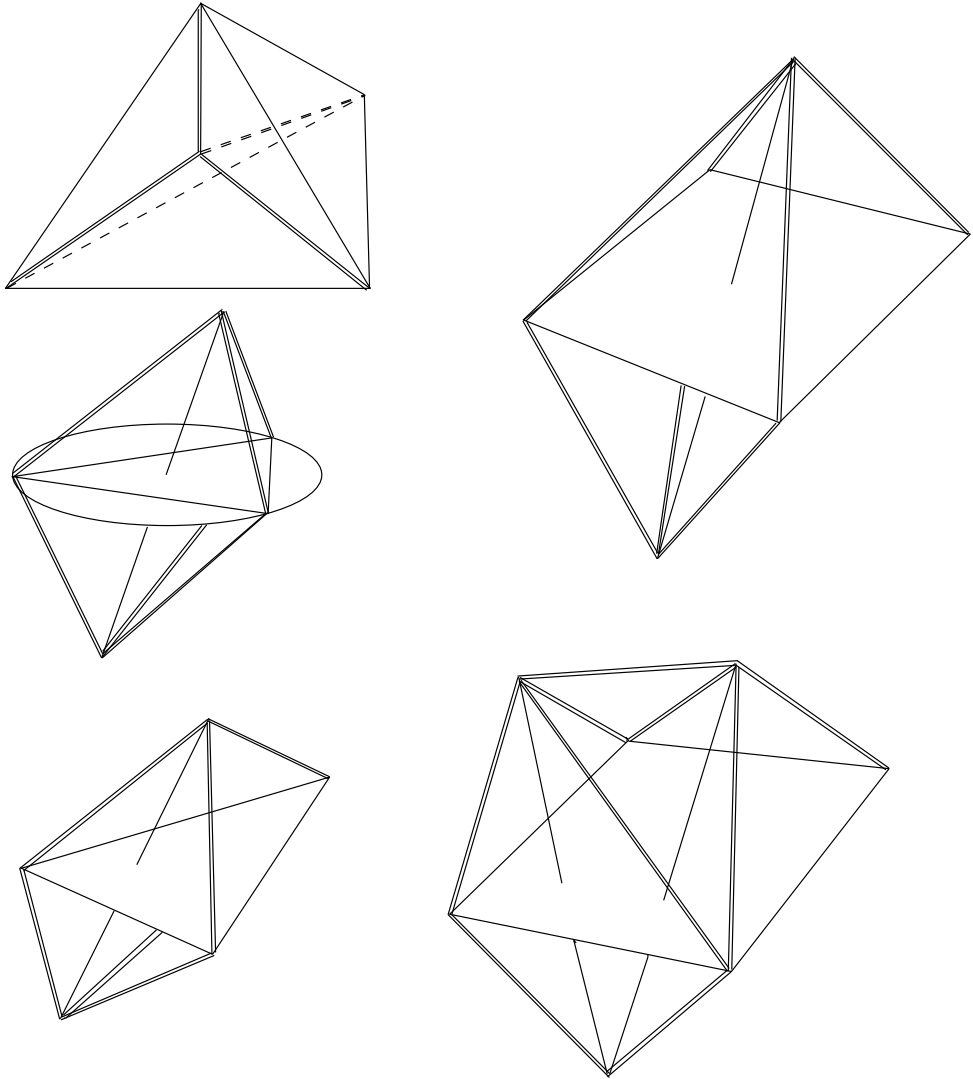}
  \caption{Some models for the sample 
	quantifier elimination problems. Struts are doubled lines.}
  \label{no-ext}
\end{figure}

(Bob Connelly has pointed out that some of these problems can be
viewed as tensegrity problems, but we do not see how to
treat them all as tensegrities, so we do not
pursue this point of view here.  Globally rigid tensegrities are
analogous to our models.  However, our models are not claimed to be
rigid.)

\begin{example}
\label{example:simplex}
In Problem~\ref{first}, let the upper bounds on the edges of the simplex
be $\sqrt8$, and let $r=2$.
We take the model to be 
a regular tetrahedron with edges marked as cables.  
Mark the edges from the circumcenter to the vertices as struts.
We apply pivots to the simplex to increase
its edges to $\sqrt8$.  
Move the interior point by a sequence of pivots
so that it has distance exactly $2$ to three of the four vertices of
the simplex.  

After these pivots are completed,  
the configuration is uniquely determined, and
a calculation with explicit coordinates shows that the configuration
does not exist, because the distance from the interior point to the fourth
vertex of the simplex is too small.
\end{example}

In these problems, when we pivot in the correct direction, the
distance constraints between points take care of themselves.  However,
some of the problems impose additional constraints.  In Problem~\ref{first},
the point is constrained to lie in the simplex.  In the other problems,
lines are required to be linked with various space polygons.
A separate verification is required to see that pivots do not violate
these additional constraints.  
{\it These separate verifications are again
quantifier elimination problems, of a smaller magnitude than the original
problem.}  

For example, in Problem~\ref{first}, we verify that
the point cannot be an interior point of a face of the simplex.
This insures that the point does not escape from the interior of the simplex
during the sequence of pivots.  The argument that there does not
exist a point in a face, under the given metric constraints, is similar
to Example~\ref{example:simplex}, but all the arguments are reduced to
two dimensions, instead of three.

The preprocessing in most other cases is similar to 
Example~\ref{example:simplex}, 
and can be reconstructed
without difficulty from the models.  The two exceptions are 
Problem~\ref{six} and
Problem~\ref{seven},
which require substantial preprocessing and a lemma to insure that the
pivots can be carried out.  (I would much prefer to have arguments based
on a pure quantifier elimination algorithm and bypass this lemma entirely,
but the current quantifier elimination algorithms do not seem up to
the task.)

\begin{lemma} Fix constants $\ell_i$, $k_i$, $\epsilon$, $h_i$, and $\ell$
subject to the constraints
	\begin{equation}
		\ell_i < \sqrt8,
			\quad k_i\in[2,2.61],
			\quad \epsilon\ge2,
			\quad h_i\in[2,\sqrt8],
			\quad\ell\in[2,\sqrt8].
	\end{equation}
Pick the following parameters in Example~\ref{seven}
$$\begin{array}{lll}
\dmin(p,q_j)=h_i & \dmin(q_1,q_2)=\epsilon & \dmin({\text{others}})=2\\
\dmax(0,p_i) = \ell_i & \dmax(0,q_j) = \ell & \dmax(p_i,p_{i+1})=k_i
\end{array}
$$
and let the other values of $\dmax$ be $+\infty$.  If a configuration
of $7$ points exists with these parameters, then a configuration also
exists with these parameters and the additional constraints
that
$$|p_i|=2,\quad |p_i-p_{i+1}| = k_i,\quad |q_j| = \ell.$$
Furthermore, the same lemma and conclusion holds in the context of
the 6-point configuration of Example~\ref{six}, if we take
$q_1=q_2=q$ and $\dmin(q_1,q_2)=0$.
\end{lemma}

\begin{proof}  This is Lemma 4.3 of \cite{i}.  Some of the constants
have been relaxed in a way that affects the proof in a very minor way.
(Two modifications must be made to the proof.  The assertion that
the circumradius of a triangle of sides $2.1,2.51,2.51$ is less than $\sqrt2$
of the original must be replaced with the
assertion that there exists an $x>2$ such that the circumradius of the
triangle of sides $x,2.61,2.61$ is less than $\sqrt2$.
Also, an instance of Problem~\ref{first} is needed, with $r=2$ and
the length-bounds of the sides of the simplex $2.61$, $2.61$, $\sqrt8$ at one
vertex, and lengths at most $\sqrt8$ at the edges opposite the edges
of length at most $2.61$.)
\end{proof}

\section{Linear Assembly Problems}
\label{linear}

In this section we define a class of nonlinear optimization problems that we
call {\it linear assembly problems}.

Assume 
given a topological space $X$, and a finite collection of topological
spaces, called {\it local domains}.  For each local domain $D$ there is a map
$\pi_D:X\to D$.  There are functions $u_i$, 
$i=1,\ldots,N$, each defined on some local domain
$D_i = \dom(u_i)$, and we let $x_i$ denote the composite 
$x_i = \pi_{D_i}\circ u_i$.

On each local domain $D$, the functions $u_i$ are related by a finite set
of nonlinear relations 
\begin{equation}\phi(u_i : \dom(u_i) = D) \ge0, \quad \phi \in \Phi_D.
	\label{phi}
\end{equation}

We use vector notation $x = (x_1,\ldots,x_N)$, with constant vectors
$c$, $b$, and matrix $A$ given.

The problem is to maximize $c\cdot x$ subject to the constraints
	\begin{equation}\label{Ax}A\, x \le b,
	\end{equation}
and to the nonlinear relations~\ref{phi}.  A problem of this form is
called a linear assembly problem.  (Intuitively, there are a number of
nonlinear objects $D$, that form the pieces of a jigsaw puzzle that
fit together according to the linear conditions~\ref{Ax}.)

\begin{example}
Assume a single local domain $D$, and let $\pi_D:X=D$ be the
identity map.  The function $f = c\cdot x $ is nonlinear.
The problem is to maximize $f$ over $D$ subject to the nonlinear relations
$\Phi_D$.  This is a general constrained nonlinear optimization problem.
\end{example}

\begin{example}  Assume that each $u_i$ has a distinct local domain $D_i = \RR$.
Let $X = \RR^N$, let $\pi_D$ be the projection onto the $i$th coordinate,
and let $x_i$ be the $i$th coordinate function on $\RR^n$.  
Assume that $\Phi_D$ is empty for each $D$.  The problem becomes
the general linear programming problem
	$$\max c\cdot x$$
such that $A x\le b$.
\end{example}

These two examples give the nonlinear and linear extremes in linear
assembly problems.
The more interesting cases are the mixed cases which combine nonlinear and
linear programming.  Example~\ref{pr:third} gives one such case.

\begin{figure}[htb]
  \centering
  \includegraphics{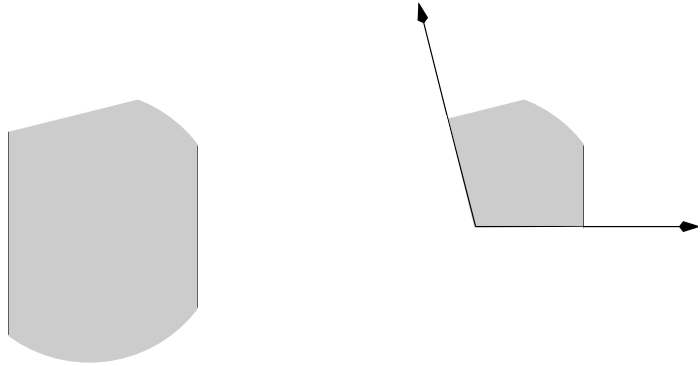}
  \caption{A truncated Voronoi cell and a subset of the cell lying in a sector}
  \label{voronoi}
\end{figure}

\begin{example} (2D Voronoi cell minimization). Take a packing of disks of
radius $1$ in the plane.  Let $\Lambda$ be the set of centers of the
disks.  Assume that the origin $0\in\Lambda$ is one of the centers.
The truncated Voronoi cell at $0$ is the set of all $x\in\RR^2$ such
that $|x|\le t$, and $x$ is closer to the origin than to any other
center in $\Lambda$.  We assume $t\in(1,\sqrt2)$.

Only the centers of distance at most $2t$ affect the shape and
area of the truncated Voronoi cell.  For each $n=0,1,2,\ldots$, we have
a topological space of all truncated Voronoi cells with $n$ nonzero
disk centers $v_i$ 
at distance at most $2t$.  Fix $n$, and let $X$ be the topological
space.

Let $D=D_i$, $i=1,\ldots,n$,  be the sectors lying between consecutive
segments $(0,v_i)$.  Each sector is characterized by its angle $\alpha$ 
and the lengths $y_a$ and $y_b$ of the two segments $(0,v_i)$,
$(0,v_j)$ between which the sector lies.  The part $A$ in $D$
of the area of the
truncated Voronoi cell is a function of the variables $\alpha$,
$y_a$, $y_b$.  A nonlinear implicit
equation $\phi=0$ relates $A$, $\alpha$, $y_a$, and $y_b$ on $D$.  
The variables $u_i$
of the linear assembly problem for the local domain $D$
are $A$, $y_a$, $y_b$, $\alpha$.

We have a linear assembly problem.  The function $c\cdot x$ is
the area of the truncated Voronoi cell, viewed as a sum of variables
$A$, for each sector $D$ (or rather, their pullbacks to $X$ under the
natural projections $X\to D$).  

The assembly constraints are all linear.
One linear relation imposes that 
the angles of the $n$ different sectors must sum to $2\pi$. Other
linear relations impose that the variable $y_a$ on $D$ equals the
variable $y_b$ on $D'$ if the two variables represent the length of
the same segment $(0,v_i)$ in $X$.
\end{example}

\subsection{Solving linear assembly problems}

In this section we describe how various linear assembly problems are
solved in the proof of the Kepler conjecture in terms sufficiently general
to apply to other linear assembly problems as well.

Let us introduce some general notation.  Let $x_D = (x_i: \dom(u_i)=D)$
be the vector of variables with local domain $D$.  Write $c\cdot x$ in
the form $\sum_D c_D\cdot x_D$ and the assembly conditions as
$$A \,x =\sum_D A_D x_D,$$
according to the local domain of the variable.

\subsubsection{Linear relaxation} 
The first general technique is {\it linear relaxation}.  
We replace the nonlinear relations $\phi(x_D)\ge0, \phi\in\Phi_D$ 
with a collection of linear inequalities that are true whenever
the constraints $\Phi_D$ are satisfied:
$A'_D x_D \le b_D$.  A linear program is obtained by replacing
the nonlinear constraints $\Phi_D$ with the linear constraints. Its
solution dominates the nonlinear optimization problem.  In this
way, the nonlinear maximization problem can be bounded from above.

Let us review some constructions that insure rigor in 
linear programming solutions.
We assume general familiarity with the basic theory 
and terminology of linear programming. 
It is well-known that the primal has a feasible solution iff the dual
is bounded.  We will formulate our linear programs in such a way
that both the primal and the dual problems
are feasible and bounded.

We use vector notation to formulate a primal problem as 
	\begin{equation}
		\max\, c\cdot x
		\label{cx}
	\end{equation}
such that $A x \le b$, where $x$ is a column vector of free variables (no
positivity constraints), $A$ is a matrix, $c$ is a row vector, and
$b$ is a column vector. 

We can insure that this primal problem is bounded by bounding each of the
variables $x_i$.  (This is easily achieved considering the geometric origins
of our problem, which provides interpretations of variables as 
particular dihedral angles,
edge lengths, and volumes.)  
We assume that these bounds form part of the constraints $A x\le b$.

The linear programs we consider
have the property that if the maximum is less than a constant $K$,
the solution does not interest us.  (For instance, in the dodecahedral
conjecture, Voronoi cell volumes are of interest
only if the volume is less than the volume of the regular dodecahedron.)
This observation allows us to replace the primal problem with 
one having an additional
variable $t$:
	\begin{equation}
		\max\, c \cdot x + K\, t
	\end{equation}
such that $A x + b t \le b$, and $0\le t \le 1$.  This 
modified primal is bounded for
the same reasons that the original primal is.  It has the feasible solution
$x=0$ and $t=1$.

\begin{lemma}
If the maximum $M$ of the original primal is greater than $K$, 
then the optimal solution of the modified primal has $t=0$, and hence its maximum is
also $M$.
\end{lemma}

\begin{proof}
Assume that $(x_0,t_0)$ gives an optimal 
solution to the modified problem
for some $1>t_0>0$, with
$c \cdot x_0 + K t_0 > M$.
Then
$(x_1,t_1) = (x_0/(1-t_0),0)$ is also a feasible solution and it beats the
optimal solution
	\begin{equation}
		c \cdot x_1 + K t_1 > c \cdot x_0 + K t_0.
	\end{equation}
This contradiction proves $t_0=0$.
\end{proof}

The output from linear program that is solved by numerical methods can
be transformed into a rigorous bound as follows. 
Based on the preceding remarks, we assume
that these linear programs are feasible and bounded.  
The dual is then also feasible and bounded.  We assume that the
numerical solutions are carried out with sufficient accuracy to insure
bounded feasible approximations to the true optima.

To explain the rigorous verification, we separate the equality constraints
from the inequality constraints, and rewrite the problem as
	\begin{equation}
		\max c \cdot x
	\end{equation}
such that $A' x = b'$, $A x \le b$, with $x$ free.
The dual problem yields a solution to
	\begin{equation}
		\min y b'  + z b,
	\end{equation}
such that $y A' + z A = c$, with $z\ge 0$ and $y$ free.
Let $(y_0,z_0)$ be a numerically obtained approximation to the dual
solution.  The vector $z_0$ will be approximately positive, and by
replacing negative coefficients by $0$, we may assume $z_0\ge0$.
Let $\delta = c - y A' - z A$ 
be the error row vector resulting from
numerical approximations.
Then for any feasible solution $x$ of the primal, we have 
	\begin{equation}
		c \cdot x 
			= (\delta + y_0 A' + z_0 A) x 
			\le \delta\cdot x + y_0 \cdot b' + z_0 \cdot b.
	\end{equation}
Using the bounds of the variables $x_i$, we bound $\delta\cdot x \le D$,
and thus obtain the rigorous upper
bound $c \cdot x \le D + y_0 \cdot b' + z_0 \cdot b$ on the primal.

\subsubsection{Implementation details}

The linear programs are solved numerically using a commercial package (CPLEX).
The input and output to these numerical programs are processed by a
custom java program, which is linked to CPLEX with a java API provided
by the software manufacturer.  
Each bound is calculated with
interval arithmetic to insure that it is reliable.
(We use a simple implementation of interval arithmetic in java based
on the math.BigDecimal implementation of arbitrary precision arithmetic.)

\subsection{Nonlinear duality}  
The second general technique is nonlinear duality.  Suppose that
we wish to show that the maximum of the primal problem~\ref{cx} 
is at most $M$.  

Let $x^* = (x^*_D)$ be a guess of the solution 
to the problem, obtained
for example, by numerical nonlinear optimization.   
We relax the nonlinear optimization by dropping from the matrix $A$
and the vector $b$ those inequalities that are not binding at $x^*$.
With this modification, we may assume that $A\,x^*=b$.  Let $m$ be
the size of the vector $b$, that is, the number of binding linear
conditions. Let $d$ be the number of local domains $D$.

We introduce
a linear dual problem with real
variables $t$, $r_\phi: \phi\in\Phi_D$, and $w\in\RR^m$.
The variables $r_\phi$ and $w$ are constrained to be non-negative.

We consider the linear problem of maximizing $t$ such that
	\begin{equation}
		M + d\, t - c\cdot x^* \ge 0
		\label{Mx}
	\end{equation}
and such that for each $x_D$ in each $D$ the linear inequality
	\begin{equation}
		c_D\cdot (x_D-x^*_D) + \sum_{\Phi_D} r_\phi \phi(x) +
					w A_D (x^*_D-x_D) + t 
			< 0
		\label{xD}
	\end{equation}
is satisfied.

There is no guarantee that a feasible solution exists to this system
of inequalities.  However, any 
feasible solution gives an upper bound $M$.
Indeed, let $x=(x_D)$ be any feasible argument to the primal, and let
$t,r_\phi,w$ be a feasible solution to the dual.
Taking the sum of the linear inequalities~\ref{xD}, 
over $D$ at $x$, we have (recall $\phi\ge0$ and $A x\le b$):
$$
\begin{array}{lll}
M &\ge M + c\cdot (x-x^*) + \sum_D\sum_{\Phi_D} r_\phi \phi(x)
	+ w A (x^*-x) + d\, t,\\
	&\ge c\cdot x + (M + d\, t - c\cdot x^*) + w (b-A x),\\
	&\ge c\cdot x.
\end{array}
$$

Since the dual problem has infinitely many constraints 
(because of constraints
for each $x\in D$), we solve the dual problem in two stages.
First, we approximate each $D$ by a finite set of test points, and solve
the finitely constrained linear programming problem for $t, r_\phi$, and $w$.

We replace $t$ with $t_0 = (-M +c\cdot x^*)/d$ (to make the constraint
\ref{Mx} bind).  It follows from
the feasibility of $t$ that $t\ge t_0$, and that 
$t_0,r_\phi,w$ is also feasible on the finitely constrained problem.
To show that $t_0,r_\phi,w$ satisfies all the inequalities~\ref{xD}
(under the substitution $t\mapsto t_0$),
we use interval arithmetic to show that each of these inequalities hold.
(To make these interval arithmetic verifications as easy as possible,
we have chosen the solution $t_0,r,w$ to make the closest inequality
hold by as large a margin $t-t_0$ as possible. This is the meaning of the
maximization over $t$
in the dual problem.)  The next section will give further details
about interval arithmetic verifications.

\subsection{Branch and bound}
The third technique is branch and bound.  When no feasible solution
is found in step (2), it may still be possible to partition $X$ into
finitely many sets $X = \coprod X_i$, on which feasible solutions to the
dual may be found.  Although this is an essential part of the solution,
the rules for branching in the Kepler conjecture follow the structure
of that problem, and we do not give a general branching algorithm.

\section{Automated Inequality Proving}
\label{inequality}

What we would like is a general, efficient algorithm for proving
inequalities of several real variables.  Each inequality $f < 0$
of a continuous function on a compact domain
can be
expressed as a maximization problem: 
	\begin{equation}
		\max f < 0.
	\end{equation}
Generally efficient algorithms are not possible because NP hard 
problems can be encoded as optimization problems of
quadratic functions
\cite{HPT}.

This section describes an inequality proving procedure that has worked
well in practice, and which could be automated to provide a
method of general interest.  This section assumes some general familiarity
with issues of floating-point and interval arithmetic, such as can be
found in \cite{A}, \cite{AH}. Our methods are similar to those in
\cite{K}.

 To prove $f<0$, it
is enough to show that the maximum of $f$ is less than $0$.  For this
reason, we use interval arithmetic to bound the maximum of 
functions.  
Through interval arithmetic, an interval $[a,b]$ containing the range of $f$
can be obtained.  By verifying that $b<0$, it follows that the range of
$f$ is negative, and hence that $f<0$.

All our functions can be built from arithmetic operations. 
(Transcendental functions are replaced with
explicit rational approximations with known error bounds.)

Often, the functions $f$ are twice continuously differentiable.
To obtain additional speed and accuracy, we use interval arithmetic
to obtain rigorous bounds on the second partial derivatives of $f$.
(We obtain formulas for the second partials through symbolic 
and automatic differentiation
of the function $f$).   With bounds on the second partials, we
obtain rigorous bounds on $f$
through its Taylor approximation.

The accuracy of the Taylor approximation improves as the domain shrinks
in size.  We chop the domain into a collection of small rectangles
and check on each rectangle whether the Taylor bound implies $f<0$.
If Taylor bound is too crude to give $f<0$, we divide it into smaller
rectangles and recompute the Taylor bounds.  By a process of adaptive
subdivision of rectangles, the inequality $f<0$ is eventually established.

Derivative information can be used to speed up the algorithm.
Taylor bounds can also be applied to the first partial derivatives of $f$.
If a partial derivative of a variable $x$ is of fixed sign on a rectangle, 
then the function is maximized along an edge $x=a$ of the rectangle.
If this edge is shared with an adjacent rectangle, the maximization of 
$f$ is pushed to an adjacent rectangle.  If this edge lies on the boundary
of the domain, the dimension of the optimization problem is reduced by one.

The method outline above works extremely well for simple functions in
a small number of variables.  The complexity grows rapidly with
the number of variables.  We are able to obtain satisfactory results
for many inequalities that depend on a single simplex $S$, that is,
functions of six variables parameterized by the edge lengths of
a simplex.

\subsection{Generative Programming}

Most of the computer code for the proof of the Kepler conjecture implements
the Taylor approximations of the nonlinear functions.
The computer code for proving $f<0$ is obtained as follows.

First, an expression for $f$ is derived.
The formulas for the first and second partial derivatives of the function
are obtained (say by a symbolic algebra system) from the expressions
for $f$. 

These symbolic expressions are then converted to an interval arithmetic
format.  In a language such as C++ with operator overloading, this
can be achieved by defining a class for intervals and overloading 
arithmetic operations so that they may be applied to instances of the
class.  In languages without operator overloading, the conversion from
the symbolic expression to computer source code is more involved.

There are other considerations to bear in mind in producing the interval
code.  
In practice, there is a substantial degradation of performance when
the rounding mode on the computer is frequently switched, and often
it is necessary to rearrange the code substantially to reduce the
number of changes in rounding mode.
Also, floating point arithmetic is not associative, so that
in order to obtain rigorous results based on interval arithmetic, great
care must be paid to the placement of parentheses.  
Another issue is the input of floating point constants.
In C++, the line of code in~\ref{x=1}
sets $x=1.0$, no matter the rounding flags.  (The constant 
is parsed at compile time and truncated to
16 digits, and there is no control over rounding modes 
until later, when the program
executes.)  The code must insure that no errors are introduced through
compiler constant truncation.
\begin{equation}
\label{x=1}
{
\text{\tt x = 1.000000000000000000001; // set x=1.0, regardless of rounding flags}
}
\end{equation}

There are many such perils in the production of reliable interval arithmetic
code.
Overall, a great
deal of effort must be expended to produce the computer code for rather
simple inequalities.  
This effort must be expended every time a new function is introduced into
an inequality.  This simple fact has kept the inequality-proving software
developed for the proof of the Kepler conjecture from having more
widespread applicability to more general inequality proving.

Figure~\ref{c++} shows a snippet of C++ code that computes
the arctangent of a linear germ of a function.

\begin{figure}[htb]
\caption{Code to calculate an interval version of the arctangent function}
\label{c++}
{

\tt

\noindent
\verb+/** + \\
\verb+ * A lineInterval is an interval version of a linear approximation to a + \\
\verb+ * function in 6 variables. + \\
\verb- * The linear approximation is +f + Df[0] x0 + Df[1] x1 + Df[2] x2 +...+ Df[5] x5. -\\
\verb+ */ + \\
\verb+class lineInterval { +\\
\verb+   public: +\\
\verb+      interval f,Df[6]; +\\
\verb+      // rest of class omitted  +\\
\verb+}; +\\
\verb+ +\\
\verb+/** +\\
\verb+ * Sample implementation of the arctangent function. +\\
\verb+ * This computes the linear approximation only.  The second derivatives +\\
\verb+ * are much more involved. +\\
\verb+ */ +\\
\verb+static lineInterval atan(lineInterval a,lineInterval b) // atan(a/b); +\\
\verb+{ +\\
\verb+   static const interval one("1"); +\\
\verb+   lineInterval temp; +\\
\verb+   temp.f = interMath::atan(a.f/b.f); // computes interval-valued arctangent +\\
\verb-   interval rden = one/(a.f*a.f+b.f*b.f); -\\
\verb!   for (int i=0;i<6;i++) temp.Df[i]= rden*(a.Df[i]*b.f-b.Df[i]*a.f); !\\
\verb+   return temp; +\\
\verb+} +\\
\verb+  +\\

}
\end{figure}

If the Kepler conjecture is eventually to be proved by generic tools, we
must find a less cumbersome way to produce the computer code.  Indeed,
a fundamental principle of software design is that there should be no
manual procedures (Pragmatic Programmer, Tip 61) \cite{HTC}.
Generative programming gives methods to automate the production
of computer code \cite{CE}.  There is nothing about the 
interval arithmetic computer code for a new function that requires
human thought or effort in an essential way.  For example, an examination
of the code for the arctangent in Figure~\ref{c++} reveals that it is
a shallow reformatting of the formula for the derivative of the
arctangent, combined with the quotient rule in calculus.
Why should the code be produced by hand, if it the process is entirely
mechanical?

A generative program could be written that takes as its input a function
and produces as output the interval arithmetic computer code
for the Taylor series bounds of that function.
The program would parse the definition of the function, generate
symbolic derivatives of the function, convert the derivative information
to computer code for calculating the derivatives, and so forth.

What advantages would this bring? First of all, it would no longer be
necessary to read 40K lines of check to check the correctness of the
proof of the Kepler conjecture. 
 It would be enough to check the code on the much smaller
generator.  Also, the same generator could be used to prove many other
inequalities. 

\subsubsection{Implementation details}
The generative program has not been written.  Some feasibility tests
have been made with javaCC for parsing and XSLT for abstract syntax
tree transformations.

\section{Plane Graph Generation}
\label{graph}
A sphere graph is a graph together with an embedding of it into
the unit sphere.  We discuss a simple sphere graph generating algorithm
in this section.

\begin{figure}[htb]
  \centering
  \includegraphics{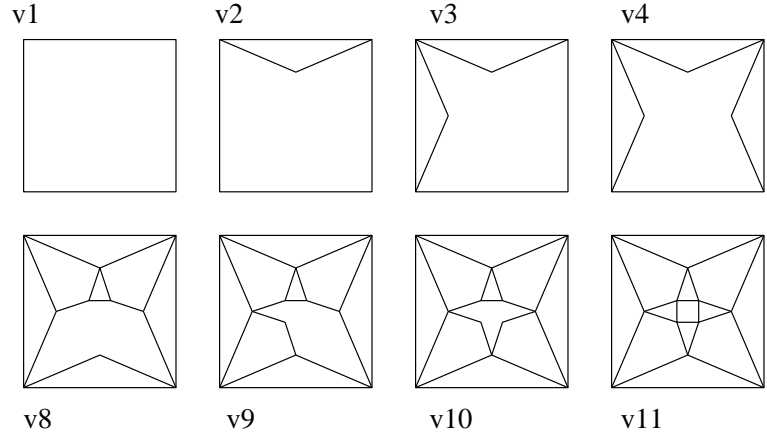}
  \caption{The first few stages and the final few stages of drawing the
	graph dual to the rhombic dodecahedron}
  \label{graph-seq}
\end{figure}

\begin{figure}[htb]
  \centering
  \includegraphics{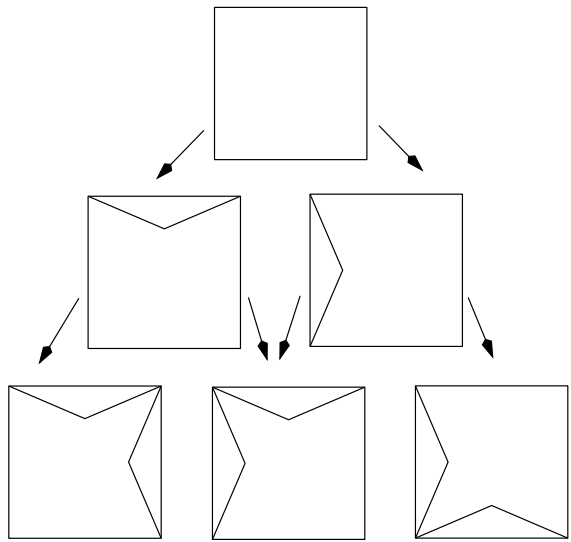}
  \caption{The same graphs can be generated by different sequences of
	adding faces}
  \label{graph2-seq}
\end{figure}

Figure~\ref{graph-seq} 
shows a sequence of sphere graphs, giving a sequence of faces
that are added to a square to produce the graph dual to the edge
graph of the rhombic dodecahedron.
We can represent
the sequence abstractly as a directed graph with vertices $v_1,\ldots,v_{11}$
with edges from $v_i$ to $v_{i+1}$.  
Figure~\ref{graph2-seq} shows that the
sequence of faces can be generated in different orders, and that all different
sequences can be represented as a directed graph 
whose root is the square $v_1$.
Call this the derivation graph.
The terminal vertices of the directed graph
represent sphere graphs isomorphic to 
$v_{11}$.

In going from a parent to a child, we always add one face.  Each
sphere graph in the derivation graph has two types of faces --
those such as the pentagon in $v_2$ that does not survive unmodified
in the terminal nodes, and those such as the triangle in $v_2$ that do.
Call these two types of faces {\it modifiable} and {\it unmodifiable}
respectively.

Let us generalize this construction to generate all the sphere graphs
that are needed for the proof of the Kepler conjecture.
Fix a natural number $N$.

Consider the set $V_0$ of nonempty sphere graphs with at most $N$ vertices, and
no loops or multiple joins.   
All faces of the sphere graph are assumed to
be polygons, and all polygons are assumed to be simple.  
We give each face one of the two attributes
{\it modifiable} or {\it unmodifiable} and call a graph with
these attributes a {\it decorated} graph.
Let $V$ be
the set of all decorated graphs of $V_0$.

Let $P$ be a polygon. We say that a simple polygon $Q$ is an 
{\it admissible refinement} of $P$ if every vertex of $Q$ is either a vertex
of $P$ or an interior point of $P$ and if $Q$ shares at least one
edge with $P$.  An admissible refinement of $P$ partitions the region
$P$ into the region $Q$ and finitely many other polygons.

Let $v$ be a decorated sphere graph in $V$.  We say that $v'$ is an
{\it admissible refinement} of  $v$ if there is a modifiable face $P$ of $v'$ 
and an admissible refinement $Q$
of $P$ such that the graph obtained by adding $Q$ is $v'$ that is compatibly
decorated.  We say that $v'$ is compatibly decorated if the unmodifiable
faces of $v'$ are $Q$ together with the unmodifiable faces of $v$.
A special case occurs, when $P=Q$, and in this case, we simply change
the attribute of $P$ to unmodifiable.

The set $V$ becomes a directed graph $\Gamma$ with an edge from each $v$ to
all of the admissible refinements of $v$.  
The root of the directed graph is an empty node.
The children of the root
are graphs consisting of a single polygon dividing
the sphere into an interior and exterior, one side modifiable and the
other not.
The terminal vertices in this directed
graph are decorated graphs, with no modifiable faces.  Thus, terminal
vertices
are in natural bijection with $V_0$.

If we take any graph in $V_0$, it can be reached from the root
as follows.  Pick a face of $V_0$ and draw its edges as
the initial polygon.  Make the interior of the
polygon modifiable.  Then continue to pick faces that have at least
one edge already drawn, and draw all remaining edges of that face,
marking the completed face as unmodifiable. This corresponds to following
a edge in the directed graph $\Gamma$.  

This gives us an algorithm to generate all sphere graphs in $V_0$: begin
with the children of the root vertex (polygons with at most $N$ vertices) and
generate all admissible refinements (that is, follow all 
possible directed edges)
until terminal vertices
are reached.  

We can improve on this algorithm by fixing for each $v\in \Gamma$
a modifiable face $P$ and an edge on that face, and then taking
only admissible refinements $Q$ of $P$ that share the given edge.
Each sphere graph in $V_0$ is still generated under this restriction.
We can also assume without loss of generality that the initial polygon
is chosen to be one with the most edges.

There is an enormous combinatorial explosion as 
all admissible refinements are generated.  As fortune has it, we are
not interested in all sphere graphs $V_0$, but rather only those that
arise as a potential counterexample to the Kepler conjecture.
Let $V_1\subset V_0$ denote this smaller set of relevant graphs.
This allows us to combine the general graph-generating algorithm with
pruning operations that keep the combinatorics from getting out of hand.

What is needed are criteria on $v\in \Gamma$ that allows us to
conclude that $v$ has no 
no admissibly refined descendents in $V_1$, that is, to conclude
there is no
directed path from $v$ to $v_1\in V_1$.  
The proof of the Kepler conjecture gives a long list of properties of
graphs in $V_1$ and this avoids the combinatorial explosion.
The implementation of the graph generator includes many other minor
tricks to keep the execution time manageable.  Pruning and the other
tricks are rather mundane, and we refer the reader to the source code
for details.

\section{Conclusion}
\label{conclusion}

We will not try to list all of the places where the algorithms of this
article would bring a simplification of the 1998 proof of the Kepler conjecture.
The list would be extensive.
To give a rough indication, we list a few places these algorithms
are relevant to the two shortest articles of the proof
\cite{FH} and \cite{iii}.

In the article~\cite{FH} alone, low-dimensional quantifier elimination problems are
the subject of Lemma 1.2, Lemma 1.3, Lemma 1.4, Lemma 1.5, Lemma 1.6,
Lemma 1.7, Lemma 1.8, Lemma 1.9, Lemma 1.11, Lemma 2.1, and Lemma 2.2.
In~\cite{iii}, an additional low-dimensional quantifier elimination problem
appears in Lemma 2.2.  In general, the parts of the 1998 proof that
rely on what that proof calls {\it geometric considerations} are
amenable to preprocessed quantifier elimination.

Linear assembly algorithms generalize the algorithm presented in
~\cite{iii}, Appendix 2.  Some examples where linear assembly would
simplify the 1998 proof are~\cite{FH} Section 4, and~\cite{iii}
Proposition 4.1, Proposition 4.2, Proposition 5.2, Proposition 5.3,
Appendix 1 (A.5 and A.7).

Interval arithmetic inequalities are used throughout the 1998 proof, in
sections such as \cite{FH} Appendix 3.13.1--4.7.5 and \cite{iii}
Section 10, Appendix 1.  In the early articles in the series,
interval arithmetic Taylor approximations are not used
\cite{i}.  As a result,
these early papers only prove very limited types of inequalities.
The entire strategy of the proof of the Kepler conjecture in \cite{i}
is shaped by these algorithmic limitations.  This profoundly affects
the structure of the optimization problem in \cite{FH}, because a
scoring function within the reach of the early primitive algorithms was
chosen, although such a function was highly suboptimal.
With improved automated inequality proving algorithms, it should be
possible to make a fresh start and devise a much more efficient scoring
function.

Graph generation is carried out in \cite{iii} Section 8.

We have not yet reached the fundamental objective
of avoiding all manual procedures; some parts of the proof
remain hand-made (even after taking account of 
the algorithms of this article).  
A fully automated proof would have to develop additional algorithms
to prove these estimates.
Nevertheless, this article brings us one step closer to
that objective.


\section*{About Authors}
Thomas C. Hales is at the Department of Mathematics, University
of Pittsburgh.

\section*{Acknowledgments}
The interval-arithmetic algorithms were developed in collaboration
with S. Ferguson.

\end{document}